\documentclass[preprint,10pt,authoryear]{elsarticle}

\usepackage{amssymb}
\usepackage{amsthm}

\newtheorem{theorem}{Theorem}
\newtheorem{lemma}[theorem]{Lemma}

\newtheorem{corollary}[theorem]{Corollary}

\theoremstyle{remark}
\newtheorem{remark}[theorem]{Remark}

\begin{document}

\begin{frontmatter}



\title{Dynamic adaptive multiple tests with finite sample FDR control}


\author{Philipp Heesen}
\ead{heesen@math.uni-duesseldorf.de}
\author{Arnold Janssen\corref{mycorrespondingauthor}} 
\ead{janssena@math.uni-duesseldorf.de}
\cortext[mycorrespondingauthor]{Corresponding author}

\address{Department of Mathematics, Heinrich-Heine-University D\"usseldorf, Universit\"atsstr. 1, 40225 D\"usseldorf, Germany}

\tnotetext[Test]{Partially supported by the Deutsche Forschungsgemeinschaft (DFG)}

\begin{abstract}
The present paper introduces new adaptive multiple tests which rely on the estimation 
of the number of true null hypotheses and which control the false discovery rate (FDR) at 
level $\alpha$ for finite sample size. 
We derive exact formulas for the FDR for a large class of adaptive multiple tests which 
apply to a new class of testing procedures. 
In the following, generalized Storey 
estimators and weighted versions are introduced and it turns out that the corresponding 
adaptive step up and step down tests control the FDR. The present results also include 
particular dynamic adaptive step wise tests which use a data dependent weighting of the 
new generalized Storey estimators. In addition, a converse of the \cite{Benjamini_Hochberg_1995} 
theorem is given. The \cite{Benjamini_Hochberg_1995} test is the only ``distribution free'' 
step up test with FDR independent of the distribution of the $p$-values of false null 
hypotheses. 
\end{abstract}


\end{frontmatter}


\section{Introduction}
\label{Section_Introduction}

Due to the multiplicity in multiple hypothesis testing, the control of a type I error rate is a serious 
problem. The familywise error rate (FWER) which is the probability of at least one false 
rejection is known to have a lack of power when the number of null hypotheses $n$ is large. Therefore, 
\cite{Benjamini_Hochberg_1995} promoted the false discovery rate (FDR) which is the expectation of the 
portion of false rejections among all rejections of a multiple test. It is often chosen as error criterion 
to control if the number of false rejections may be a reasonable portion of all rejections. This is, for 
example, the case for genome wide association studies when the significant genes will be judged again by a 
follow-up study and for exploratory data analysis. 

A popular multiple test with proven FDR control is the linear step up test of \cite{Benjamini_Hochberg_1995} 
which relies on the Simes test for the global null hypothesis, cf. \cite{Simes_1986}. It's FDR control 
has been proven under an independence assumption, cf. \cite{Benjamini_Hochberg_1995}, and when the test 
statistics are positively regression dependent on the subset of true null hypotheses, cf. \cite{Benjamini_Yekutieli_2001}. 
It is well known that the predetermined FDR level is not exhausted by the linear step up test if there 
is at least one false null hypothesis. But the exhaustion of the predetermined FDR level is essential 
for the maximization of the power of the multiple tests. Therefore, adaptive step up tests have been 
proposed which include an estimation of the number of true null hypotheses $n_0$. A popular adaptive 
step up test is the one of \cite{Storey2005} which is related to the results of \cite{Schweder_Spjotvoll_1982}. 
The adaptive test is based on the so-called Storey estimator 
\begin{equation}
 \hat n_0 = n \frac{1-\hat F_n(\lambda)+\frac{1}{n}}{1-\lambda}, 
\end{equation}
where $0<\lambda<1$ is some tuning parameter and $\hat F_n$ is the empirical cumulative distribution 
function (ecdf) of the $p$-values $p_1,\ldots,p_n$ corresponding to $n$ null hypotheses $\mathcal{H}_i$, 
$i=1,\ldots,n$. A reasonable adjustment of $\lambda$ has been discussed by 
\cite{Storey_Tibshirani_2003}, 
\cite{Storey2005}, 
\cite{Langaas_Lindqvist_Ferkingstad_2005} 
and \cite{Liang_Nettleton_2012}. 
Further estimators have, for instance, been considered by 
\cite{Benjamini_Hochberg_2000}, 
\cite{Benjamini_Krieger_Yekutieli_2006}, 
\cite{Blanchard_Roquain_2009}, 
\cite{Celisse_Robin_2010}, 
\cite{Chen_Doerge_2012}, 
\cite{Meinshausen_Rice_2006} 
and \cite{Zeisel_Zuk_Domany_2011}. 
It is well known that the adaptive step up test corresponding to the Storey estimator leads to finite sample 
FDR control under an independence assumption, cf. \cite{Storey2005} and \cite{Liang_Nettleton_2012}. However, 
the Storey estimator considers the ecdf $\hat F_n$ only at a fixed point $\lambda$ and the data dependent 
adjustments of $\lambda$ usually do not lead to finite sample FDR control of the adaptive SU test. In the 
following, we develop new estimators which include the information of the ecdf $\hat F_n$ at more than one 
point $\lambda$. The development of the new estimators is always done in view of finite sample FDR control 
and is based on a new FDR control condition which has some similarities to the one of \cite{Sarkar_2008}. 
By including $\hat F_n$ at several so-called inspection points the resulting estimators get more robust. 
As we will see, the choice of the estimator may even be data dependent while still controlling the FDR at 
finite sample size. We will refer to the corresponding adaptive step up test with data dependent estimator 
as dynamic adaptive step up tests. Again, our results hold under an independence assumption.

Adaptive tests with control of the FWER have been considered by \cite{Finner_Gontscharuk} and 
\cite{Finner_Sarkar_Guo_2012}.

The present paper is organized as follows. Section \ref{Section_Preliminaries} sets up the basic model 
assumptions and notation. In Section \ref{Section_FDR_control} we develop a new FDR control condition and 
introduce the class of generalized Storey estimators. 
Section \ref{Section_Stationary_approach} and \ref{Section_Dynamic_approach} are devoted to the development 
of new estimators which still control the predetermined FDR level. Section \ref{Section_Adaptive_SD_tests} 
briefly extends our results to adaptive step down tests. Additionally, in Section \ref{Section_Converse} we 
give a converse result of \cite{Benjamini_Hochberg_1995}. 
Finally, all proofs and some technical results are outlined in Section \ref{Section_Appendix}.

\section{Preliminaries}
\label{Section_Preliminaries}

In this section we set up our model assumptions and recall the definitions of step up (SU) tests and the 
false discovery rate (FDR). 

First, let us introduce the {\bf Basic Independence (BI) Model} which is given as follows. Consider the 
null hypotheses $\mathcal{H}_i$, $i=1,\ldots,n$, which may randomly be true or false. Furthermore, let 
\begin{equation}\label{DefinitionModel001}
 (H_i,\xi_i)_{i \leq n} : \Omega \longrightarrow (\{0,1\} \times [0,1])^n 
\end{equation}
be an arbitrary multivariate random variable and let $U_1,\ldots,U_n$ be independent, uniformly distributed 
on $(0,1)$ and independent of $(H_i,\xi_i)_{i \leq n}$. The null hypothesis $\mathcal{H}_i$ is true if $H_i=0$ 
holds and false if $H_i=1$ holds. For convenience we will also talk about true and false $p$-values instead 
of $p$-values of true and false null hypotheses and true $p$-values are identified with the null hypotheses. 
Then $(U_i)_{i \leq n}$ denotes the vector of possible true 
$p$-values and $(\xi_i)_{i \leq n}$ the vector of possible false $p$-values. Depending on the status $H_i$, 
the $p$-value $p_i$ of the $i$-th null hypothesis $\mathcal{H}_i$ is given by $U_i$ or $\xi_i$, i.e. we have 
\begin{equation}\label{DefinitionModel002}
 p_i = (1-H_i) \cdot U_i + H_i \cdot \xi_i, \quad i=1,\ldots,n. 
\end{equation}
Moreover, let us denote the random number of true null hypotheses by 
\begin{equation}\label{DefinitionModel003}
 N_0 = \sum_{i=1}^n (1-H_i). 
\end{equation}
To avoid trivial cases let $E(N_0)$ be always positive. All subsequently considered multiple tests basically rely 
on the ecdf $\hat F_n(t)= \frac{1}{n} \sum_{i=1}^n 1\{p_i \leq t\}$, $t \in [0,1]$, and the order statistics of the 
$p_i$'s, respectively. The order statistics are denoted by $p_{1:n} \leq \ldots \leq p_{n:n}$.

The BI Model includes a unifying approach of some standard models which are used for the analysis of the FDR of 
multiple tests. 
In many cases models with deterministic but unknown indicators $(H_1,\ldots,H_n)$ are studied which code the true 
null hypotheses. Then $N_0$ is fixed and it is frequently estimated. Moreover, for the analysis of adaptive SU tests 
it is mostly assumed that the $p$-values of true null hypotheses are independent, uniformly distributed on $(0,1)$ 
and independent of the $p$-values of false null hypotheses. 
Another standard model with random $N_0$ is Efron's two groups mixture model, cf. \cite{Efron_2001}, which 
has originally been formulated for test statistics. In terms of $p$-values, it is given by the following 
submodel of the BI Model, where the random vectors $(H_i)_{i \leq n},(U_i)_{i \leq n}$ and $(\xi_i)_{i \leq n}$ 
are independent, $H_1,\ldots,H_n$ are i.i.d. Bernoulli distributed and where $\xi_1,\ldots,\xi_n$ are i.i.d. 
according to some alternative distribution function (df) $F_1$. The BI Model also covers different distributions 
and dependence structures of (\ref{DefinitionModel001}). 
 
\bigskip 
Multiple tests given by a vector of critical values have often the following structure. 
Let $\hat \alpha_{i:n}$, $i=0,\ldots,n$, be data dependent critical values with 
$0 = \hat \alpha_{0:n} < \hat \alpha_{1:n} \leq \ldots \leq \hat \alpha_{n:n} < 1$. Then the adaptive step up 
(SU) test rejects all null hypotheses corresponding to the $p$-values $p_{i:n}$ which fulfill 
$p_{i:n} \leq \hat\alpha_{R:n}$ with 
\begin{equation}
 R = \max\{i : p_{i:n} \leq \hat\alpha_{i:n}\} 
\end{equation}
and the convention $\max \emptyset =0$. If the critical values are deterministic, then this 
multiple test is just called step up (SU) test. Note that $R$ denotes the number of rejections 
of the adaptive SU test. Furthermore, let 
\begin{equation}
 V = \sum_{i=1}^n 1\{p_i \leq \hat \alpha_{R:n}, H_i=0\}
\end{equation}
be the number of falsely rejected true null hypotheses. Then the false discovery rate (FDR) of 
the multiple test is given by $FDR = E\left( \frac{V}{R\vee1} \right)$, where $R\vee1=\max(R,1)$. 

The well known SU test of Benjamini and Hochberg (BH test) is based on linear critical values 
$\alpha_{i:n}=\frac{i}{n}\alpha$, $i=1,\ldots,n$. Let us shortly assume the BI Model with 
fixed $N_0=n_0$. Several authors showed that the BH test controls the FDR at fixed level 
$0<\alpha<1$, i.e. we have $FDR\leq\alpha$. To be more precise, \cite{Benjamini_Hochberg_1995} 
showed that $FDR\leq\frac{n_0}{n}\alpha$ holds and \cite{Finner_Roters_2001} and \cite{Sarkar_2002} 
proved the equation $FDR = \frac{n_0}{n}\alpha$. Under the BI Model, it is easily seen that 
$FDR=\frac{E(N_0)}{n}\alpha$ holds for the BH test. This lack of exhaustion of the predetermined 
FDR level $\alpha$ motivates the use of adaptive SU tests which include an estimator $\hat n_0$ 
of $n_0$. These tests basically use the data dependent level $\frac{n}{\hat n_0}\alpha$ instead 
of $\alpha$ which yields the heuristic 
\begin{equation}
 FDR \approx \frac{n_0}{n} \cdot \frac{n}{\hat n_0} \alpha \approx \alpha. 
\end{equation}
The adaptive SU test of \cite{Storey2005} is based on the critical values 
\begin{equation}
 \hat \alpha_{i:n} = \left( \frac{i}{\hat n_0(\lambda)} \alpha \right) \wedge \lambda
\end{equation}
with estimator 
\begin{equation}\label{StoreyEstimator}
 \hat n_0(\lambda) = n\frac{1-\hat F_n(\lambda)+\frac{1}{n}}{1-\lambda} 
\end{equation}
for $n_0$, where $0 < \lambda < 1$ is some tuning parameter and $a \wedge b = \min(a,b)$. 
Under the BI Model with fixed $N_0=n_0$ it has been shown that this adaptive SU test still 
controls the FDR, cf. \cite{Storey2005} and \cite{Liang_Nettleton_2012}. The FDR control of 
the adaptive SU test of \cite{Storey2005} carries over to the BI Model. Here, $\hat n_0(\lambda)$ 
may be regarded as estimator for $N_0$ and $E(N_0)$, respectively. A similar motivation and 
heuristic holds again.

\section{FDR control of adaptive SU tests}
\label{Section_FDR_control}

\cite{Storey2005} mainly consider the so-called Storey estimator (\ref{StoreyEstimator}) for adaptive 
SU tests with finite sample FDR control. If one is interested in further or entire classes of estimators, 
then it is reasonable to develop a central condition for the estimators which ensures FDR control of the 
corresponding adaptive SU tests. \cite{Benjamini_Krieger_Yekutieli_2006}, \cite{Sarkar_2008} and 
\cite{Zeisel_Zuk_Domany_2011} already introduced such conditions for several classes of estimators for instance. 
In this section we derive an exact expression for the FDR and we impose a general condition which works for a 
new class of generalized Storey estimators 
and leads to finite sample FDR control of the corresponding adaptive SU tests. Moreover we give a detailed 
comparison between our condition and the one of \cite{Sarkar_2008}.

Throughout, the range of $p$-values $[0,1]$ is divided in a rejection region $[0,\lambda]$ and an 
estimation region $[\lambda,1]$ for the unknown quantity $N_0$. The subsequent adaptive SU tests 
are based on the critical values 
\begin{equation}\label{aSUcriticalvalues}
 \hat \alpha_{i:n} = \left(\frac{i}{\hat n_0} \alpha\right) \wedge \lambda, \quad i=1,\ldots,n, 
\end{equation}
where $0<\lambda<1$ is some tuning parameter and predetermined level $0<\alpha<\lambda$. The critical 
values rely on an estimator 
\begin{equation}\label{aSUestimator}
 \hat n_0 = \hat n_0((\hat F_n(t))_{t\geq\lambda}) > 0
\end{equation}
of $N_0$ which is given by a measurable function of $(\hat F_n(t))_{t\geq\lambda}$. 
Moreover, we introduce the quantities 
\begin{equation}\label{DefinitionRandV}
 R(t) = n \hat F_n(t) \quad \mbox{and} \quad V(t) = \sum_{i=1}^n 1\{p_i \leq t, H_i=0\}, 
\end{equation}
$t\in[0,1]$, as the number of all $p$-values less or equal to the fixed threshold $t$ and the number 
of true $p$-values less or equal to $t$, respectively.

\begin{theorem}\label{TheoremCondition}
 Consider the BI Model and the adaptive SU test with critical values (\ref{aSUcriticalvalues}) and 
estimator (\ref{aSUestimator}). Then  
\begin{equation}\label{TheoremCondition001}
 FDR = \frac{\alpha}{\lambda} \cdot E\left( V(\lambda) \cdot \min\left\{\frac{1}{\hat n_0},\frac{1}{R(\lambda)\alpha}\right\} \right) 
\end{equation}
holds. Consequently, the condition  
\begin{equation}\label{TheoremCondition002}
 E\left(\frac{V(\lambda)}{\hat n_0}\right) \leq \lambda
\end{equation}
implies FDR control at level $\alpha$, that means $FDR \leq \alpha$. 
\end{theorem}

To our best knowledge exact FDR formulas for adaptive tests like (\ref{TheoremCondition001}) have not be derived 
earlier under the present generality. Theorem \ref{TheoremCondition} is the key for the treatment of adaptive SU 
tests with generalized Storey estimators (\ref{GeneralizedStoreyEstimator}) which have already been mentioned in 
passing by \cite{Liang_Nettleton_2012} 
in another context without focus on finite sample FDR control. Let $0< \lambda \leq \lambda_1 < \gamma_1 \leq 1$ be tuning 
parameters. Then the generalized Storey estimator is given by 
\begin{equation}\label{GeneralizedStoreyEstimator}
 \hat n_0(\lambda_1,\gamma_1) = n \frac{ \hat F_n(\gamma_1) - \hat F_n(\lambda_1) + \frac{1}{n} }{\gamma_1-\lambda_1}  
\end{equation}
and the choice of $\gamma_1=1$ just leads to the Storey estimator (\ref{StoreyEstimator}). 

\begin{theorem}\label{TheoremGeneralizedStoreyEstimator}
 Consider the BI Model and let $0< \lambda \leq \lambda_1 < \gamma_1 \leq 1$. Then 
\begin{equation}
 E\left( \frac{V(\lambda)}{\hat n_0(\lambda_1,\gamma_1)} \Big| (\hat F_n(t))_{t \geq \gamma_1} \right) \leq \lambda 
\end{equation}
and (\ref{TheoremCondition002}) hold and the adaptive SU test with critical values (\ref{aSUcriticalvalues}) and 
generalized Storey estimator (\ref{GeneralizedStoreyEstimator}) has $FDR\leq\alpha$. 
\end{theorem}

\begin{remark}
 Theorem \ref{TheoremGeneralizedStoreyEstimator} remains true when the estimator $\hat n_0(\lambda_1,\gamma_1)$ 
is replaced by $\hat n_0(\lambda_1,\gamma_1) \cdot(1-(\frac{\lambda_1}{\gamma_1})^{R(\gamma_1)\vee1})$. The 
effect of the additional factor is vanishingly small for large $n$. Hence, we omit this factor. 
\end{remark}

The conditions of \cite{Benjamini_Krieger_Yekutieli_2006}, \cite{Sarkar_2008} and \cite{Zeisel_Zuk_Domany_2011}, 
which ensure finite sample FDR control, resemble each other. Thus, we restrict ourselves to the comparison to the methods 
of \cite{Sarkar_2008}. First, \cite{Sarkar_2008} does not take the minimum with $\lambda$ in (\ref{aSUcriticalvalues}). 
The adaptive SU tests may reject $p$-values on the entire interval $[0,1]$. However, mostly only hypotheses with 
$p$-values smaller than $\lambda$ should be rejected for practical reasons. 
The tuning parameter $\lambda$ is often chosen close to $0.5$, 
cf. \cite{Storey_Tibshirani_2003}. Furthermore, the estimators may use the information of the complete ecdf 
$\hat F_n$. But since $\hat F_n(t)$ is often biased for small $t$ by the false $p$-values, the most useful 
information should be provided by $\hat F_n(t)$ for $t\geq\lambda$. 

The condition of \cite{Sarkar_2008} only applies to estimators $\hat n_0$ which are 
non-decreasing functions in each $p$-value $p_i$. Observe that the generalized Storey 
estimators (\ref{GeneralizedStoreyEstimator}) with $\gamma_1<1$ are not non-decreasing, 
but they are very useful as we will see below. The 
methods of \cite{Benjamini_Krieger_Yekutieli_2006} and \cite{Zeisel_Zuk_Domany_2011} also 
rely on estimators which are non-decreasing functions in each $p$-value $p_i$ and are thus 
not applicable. In Section \ref{Section_Stationary_approach} and \ref{Section_Dynamic_approach} 
we will develop new weighted versions of the generalized Storey estimators and show that the 
corresponding adaptive SU tests still control the FDR. These weighted versions include the 
information of the ecdf $\hat F_n$ at several inspection points and may show a more robust behavior 
than the Storey estimator. As we will see, the weighting may even be data dependent and the not 
non-decreasing generalized Storey estimators just offer the latitude which is needed for the 
FDR control of the corresponding dynamic adaptive SU tests.

\bigskip 
In the next step we offer another method for the verification of the inequality (\ref{TheoremCondition002}) 
which is comparable with earlier work. As Lemma \ref{LemmaComparisonCondition} will show, the condition of 
\cite{Sarkar_2008} and (\ref{TheoremCondition002}) are based on the same term which are applied to different 
classes of estimators. Therefore, let $p^{(i)}=(p_1,\ldots,p_{i-1},0,p_{i+1},\ldots,p_n)$ be the vector of 
$p$-values, where the $i$-th $p$-value is set to zero and 
\begin{equation}\label{EstimatorDecreased}
 \hat n_0^{(i)} = \hat n_0(p^{(i)})
\end{equation}
be the estimator $\hat n_0$ based on $p^{(i)}$. In actual fact, we introduced $\hat n_0$ in (\ref{aSUestimator}) 
as function of a part of the ecdf $\hat F_n$. This still applies, but for notational convenience let us write 
(\ref{EstimatorDecreased}). 

\begin{lemma}\label{LemmaComparisonCondition}
 Under the assumptions of Theorem \ref{TheoremCondition}, 
\begin{equation}\label{LemmaComparisonCondition001}
 \frac{1}{\lambda} \cdot E\left( \frac{V(\lambda)}{\hat n_0} \right) = E\left( \sum_{i: H_i=0} \frac{1}{\hat n_0^{(i)}} \right). 
\end{equation}
\end{lemma}

In his set up, \cite{Sarkar_2008} showed that the condition $E( \sum_{i: H_i=0} 1 / \hat n_0^{(i)} ) \leq \alpha$ 
is sufficient for finite sample FDR control. Hence, under the assumptions of Theorem \ref{TheoremCondition}, this 
condition is just equivalent to (\ref{TheoremCondition002}). But note that Theorem \ref{TheoremCondition} applies 
to a different class of estimators. Moreover, 
observe that the left hand side of (\ref{TheoremCondition002}) factors the distribution of true $p$-values into 
the condition. Here, the true $p$-values are uniformly distributed on the unit interval. Hence, this does not matter. 
Under certain regularity assumptions, (\ref{TheoremCondition002}) also provides FDR control if the true $p$-values 
are stochastically larger than the uniform distribution on the unit interval and the condition turns out to be sharper, 
cf. \cite{Heesen_2014}. Such distributions of true $p$-values widely occur in one-sided hypothesis testing problems, 
for instance.

\section{Stationary approach}
\label{Section_Stationary_approach}

For each family of estimators (including the Storey and generalized Storey estimators) it is not 
yet clear for which estimator an adequate bias and variance may be achieved since the distribution 
of false $p$-values is unknown. Therefore, a weighted estimator may perform better. 


\begin{corollary}\label{CorollaryStationaryApproach}
Consider the BI Model, let $\beta_1,\ldots,\beta_k$ be fixed positive weights with $\sum_{i=1}^k \beta_i =1$ 
and let $\hat n_{0,i}$ be estimators which fulfill (\ref{aSUestimator}) and 
\begin{equation}
 E\left( \frac{V(\lambda)}{\hat n_{0,i}} \right) \leq \lambda \quad i=1,\ldots,k. 
\end{equation}
Then the weighted estimator
\begin{equation}\label{WeightedEstimatorFixed}
 \hat n_0 = \sum_{i=1}^k \beta_i \hat n_{0,i}
\end{equation}
satisfies (\ref{TheoremCondition002}) and the adaptive SU test with critical values (\ref{aSUcriticalvalues}) 
and estimator (\ref{WeightedEstimatorFixed}) has $FDR\leq\alpha$. 
\end{corollary}

\begin{remark}
Let $(\hat n_{0,x})_{x \in X}$ be a family of estimators (\ref{aSUestimator}). Then 
Corollary \ref{CorollaryStationaryApproach} may easily be generalized by considering 
estimators $\hat n_0 = \int \hat n_{0,x} d\nu(x)$, where $\nu$ is a probability measure 
on $X$. 
These results follow from Jensen's inequality. 
\end{remark}

By Theorem \ref{TheoremGeneralizedStoreyEstimator} it is clear that Corollary \ref{CorollaryStationaryApproach} 
is applicable for the Storey estimators $\hat n_0(\lambda_i)$ with $0<\lambda\leq\lambda_1<\ldots,\lambda_k<1$. As 
already mentioned, the Storey estimator $\hat n_0(\lambda)$ considers the ecdf $\hat F_n$ only at a fixed point 
$\lambda$ and is hence particularly sensitive. 

\bigskip
{\bf\noindent A practical guide for the weighting of the Storey estimators (\ref{StoreyEstimator}).} 
An ad hoc method for the choice of fixed weights of the weighted Storey estimator 
\begin{equation}\label{WeightedStoreyEstimatorFixed}
\hat n_0 = \sum_{i=1}^k \beta_i \hat n_0(\lambda_i)
\end{equation}
may be that the conditional variance $Var(\beta_i \hat n_0(\lambda_i)|(H_i,\xi_i)_{i \leq n})$ is constant 
for all $i=1,\ldots,k$. Then the variance formula for binomials leads to 
\begin{equation}
 \beta_i = \frac{\sqrt{\frac{1}{\lambda_i}-1}}{\sum_{j=1}^k \sqrt{\frac{1}{\lambda_j}-1}}, \quad i=1,\ldots,k. 
\end{equation}

\begin{figure}[ht]
  \centering
\includegraphics[scale=0.45]{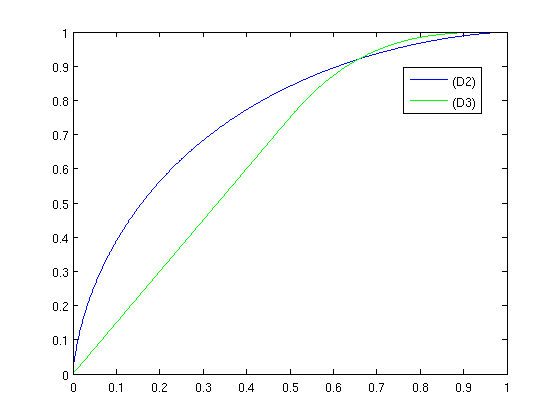}
  \caption{Df of false $p$-values in situation (D2) and (D3).}
   \label{Figure:002}
\end{figure}

We give a small simulation study which compares the false discovery rates of the adaptive SU tests with 
Storey estimator (\ref{StoreyEstimator}) and weighted Storey estimator (\ref{WeightedStoreyEstimatorFixed}) 
by a Monte-Carlo simulation with $10.000$ iterations. Therefore, consider the BI Model, let $n=1000$ and 
$N_0=600$ be fixed. In addition, let $\lambda=0.5$ be the tuning parameter of the Storey estimator, 
$\lambda_1=0.5, \lambda_2=0.6, \lambda_3=0.7$ be the tuning parameters of the weighted Storey estimator and 
$\alpha=0.05$ be the predetermined FDR level. Then the practical guide advises to set 
$\beta_1=0.4, \beta_2=0.33$ and $\beta_3=0.27$. The $n-N_0=400$ false $p$-values $\xi_i$ are i.i.d. and 
\begin{itemize}
 \item[(D1)] distributed according to the Dirac distribution $\delta_0$ with point mass at $0$, 
 \item[(D2)] given by $\xi_i=1-\Phi(X_i+1)$, where $X_i \sim N(0,1)$ with df $\Phi$, or 
 \item[(D3)] distributed according to the df 
\begin{equation}
  F_1(t) = 
 \frac{3}{2}t \cdot 1\left\{t \leq \frac{1}{2}\right\} + (1-2(1-t)^3) \cdot 1\left\{t >\frac{1}{2}\right\}.
\end{equation}
\end{itemize}
The df of false $p$-values in situation (D2) and (D3) are plotted in Figure \ref{Figure:002}. 
Note that (D1) may be replaced by 
\begin{itemize}
 \item[(D4)]distributed according to an df $F_1$ with $F_1(0.5)=1$
\end{itemize}
without changing the results of the simulation, because $0.5=\lambda=\lambda_1$ holds and given 
$(\hat F_n(t))_{t\geq\lambda}$ the adaptive SU test is just a BH test whose FDR only depends on 
the number of true null hypotheses. The results of the simulation are given in Table \ref{Table:001}. 
\begin{table}[ht]
  \centering
\begin{tabular}{c|ccc}
	$\hat n_0$	& (D1)  & (D2)  & (D3)  \\ \hline
$\hat n_0(0.5)$ & 0.0501 & 0.0392 & 0.0354 \\
$\sum_{i=1}^3 \beta_i \hat n_0(\lambda_i)$ & 0.0499 & 0.0432 & 0.0393 \\
$\sum_{i=1}^6 \hat\beta_i \hat n_0(\lambda_{i-1},\lambda_i)$ & 0.0491 & 0.0437 & 0.0434
\end{tabular}
  \caption{FDR of the adaptive SU tests with Storey estimator (\ref{StoreyEstimator}), 
weighted Storey estimator (\ref{WeightedStoreyEstimatorFixed}) and dynamic estimator 
(\ref{EsimatorDynamic}) which is introduced in Section \ref{Section_Dynamic_approach}.} 
   \label{Table:001}
\end{table}
The FDR of the classical BH test is here just $\frac{N_0}{n}\alpha = 0.03$ which is enlarged by the adaptive tests. 
(D1) corresponds to a situation with maximal signal strength of false null hypotheses and (D4) includes 
many situations with large signal strength. In these situations, both estimators essentially lead to the 
same predetermined FDR level $\alpha$. (D2) and (D3) are situations with moderate signal strength. The approach 
of (D2) is often used to model true and false $p$-values and (D3) corresponds to a non-parametric approach. Here, 
the weighted Storey estimator is significantly superior to the Storey estimator. However, under the non-parametric 
approach it is easy to find distributions of false $p$-values, where the weighted Storey estimator is even more 
superior.

\section{Dynamic approach}
\label{Section_Dynamic_approach}

Let $0<\lambda = \lambda_0 < \lambda_1 < \ldots < \lambda_k=1$ and $\beta_i= \frac{\lambda_i-\lambda_{i-1}}{1-\lambda}$, 
$i=1,\ldots,k$. Then the primarily proposed estimator of \cite{Storey2002,Storey2003} may be decomposed by 
\begin{equation}\label{Pre-weighting}
 n\frac{1-\hat F_n(\lambda)}{1-\lambda} 
= n\sum_{i=1}^k \beta_i \frac{\hat F_n(\lambda_i) - \hat F_n(\lambda_{i-1})}{\lambda_i-\lambda_{i-1}}, 
\end{equation}
where the additional term $\frac{1}{n}$ of the Storey and generalized Storey estimators is omitted. Similar 
as in Section \ref{Section_Stationary_approach}, the individual estimators $\hat n_0(\lambda_{i-1},\lambda_i)$ 
and $n\frac{\hat F_n(\lambda_i)-\hat F_n(\lambda_{i-1})}{\lambda_i-\lambda_{i-1}}$ may have an unknown bias 
depending on the distribution of false $p$-values. Furthermore, a small estimate of $N_0$ is 
desirable for a large power of the adaptive SU test. Therefore, a look at the data would be very helpful to 
adjust the the weights $\beta_i$ in this manner. Of course, FDR control may not be achieved by an arbitrary 
adjustment. But it turns out that a special data dependent weighting which uses the tail information of the 
ecdf $\hat F_n$ still leads to finite sample FDR control of the corresponding dynamic adaptive SU test and is not limited 
to the generalized Storey estimators. 

\begin{theorem}\label{TheoremDynamicApproach}
 Consider the BI Model, let $0<\lambda<\gamma_1\leq \ldots\leq \gamma_k\leq 1$ and 
$\mathcal{F}_{\gamma_i}=\sigma((\hat F_n(t))_{t\geq\gamma_i})$ be the generated 
$\sigma$-algebra, $i=1,\ldots,k$. Let $\hat \beta_i$ be non-negative, data dependent 
and $\mathcal{F}_{\gamma_i}$-measurable weights with $\sum_{i=1}^k \hat \beta_i=1$. 
Moreover, let $\widetilde n_0(\gamma_i)$, $i=1,\ldots,k$, be $\mathcal{F}_\lambda$-measurable 
 estimators of the form 
(\ref{aSUestimator}) which satisfy 
\begin{equation}\label{TheoremDynamicApproachCondition}
 E\left( \frac{V(\lambda)}{\widetilde n_0(\gamma_i)} \Big| \mathcal{F}_{\gamma_i} \right) \leq \lambda, 
\end{equation}
including the generalized Storey estimators $\hat n_0(\lambda_i,\gamma_i)$ with 
$\lambda\leq\lambda_i<\gamma_i$. 
Then the estimator 
\begin{equation}\label{DataDependentWeightedEstimatorGeneral}
 \hat n_0 = \sum_{i=1}^k \hat \beta_i \widetilde n_0(\gamma_i) 
\end{equation}
satisfies (\ref{TheoremCondition002}) and the adaptive SU test with critical values 
(\ref{aSUcriticalvalues}) and estimator (\ref{DataDependentWeightedEstimatorGeneral}) 
has $FDR\leq\alpha$.
\end{theorem}

\begin{figure}[ht]
  \centering
\includegraphics[scale=0.5]{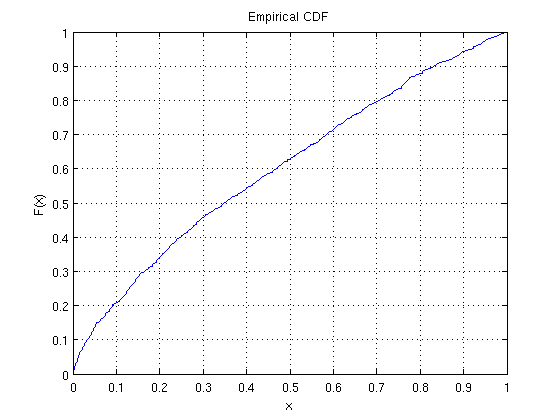}
  \caption{Likely ecdf of $p_1,\ldots,p_n$ in situation (D2).}
   \label{Figure:001}
\end{figure}
Below a practical adjustment of data driven weights $\hat \beta_i$ is motivated for the generalized 
Storey estimators (\ref{GeneralizedStoreyEstimator}). 
As illustration consider a simulation of situation (D2) of Section \ref{Section_Stationary_approach}, 
where the ecdf $\hat F_n$ is presented in Figure \ref{Figure:001}. We are looking for an interval 
$[1-\delta,1]$, where $\hat F_n$ is most informative for the estimation of $N_0$. On this part the 
generalized Storey estimator $\hat n_0(\widetilde \lambda, \widetilde \gamma)$, 
$1-\delta \leq \widetilde \lambda < \widetilde \gamma<1$, may perform good. Since $\delta$ is unknown 
the interval $[\lambda,1]$ is first divided into small pieces $[\lambda_i,\lambda_{i+1}]$. This 
observation motivates the following backwards data based dynamic adjustment of $\hat \beta_i$  
which are related to the deterministic weights $\beta_i = \frac{\lambda_i - \lambda_{i-1}}{1-\lambda}$ 
used in (\ref{Pre-weighting}). We are going to fill the last row of Table \ref{Table:001} for the 
following concrete adaptive weights $\hat \beta_i$ of that example. 
Fix in advance $\lambda=\lambda_0=0.5, \lambda_1=0.6, \lambda_2=0.7, \lambda_3=0.8, \lambda_4=0.9, \lambda_5=0.95, \lambda_6=1$ 
and let $\epsilon>0$ be a tuning parameter. The estimator 
\begin{equation}\label{EsimatorDynamic}
 \hat n_0 = \sum_{i=1}^6 \hat\beta_i \hat n_0(\lambda_{i-1},\lambda_i) 
\end{equation}
is established as follows. Based on the pre-weighting in (\ref{Pre-weighting}) define  
$$ \hat \beta_6 = \beta_6 = \frac{\lambda_6-\lambda_5}{1-\lambda}, \quad \mbox{and} \quad 
\hat\beta_5 = \beta_5 = \frac{\lambda_5-\lambda_4}{1-\lambda}.
$$
The other weights are set as follows. \\
Case 1. Assume that there is an index $1\leq i\leq4$ with 
\begin{eqnarray}
& \hat n_0(\lambda_i,\lambda_{i+1}) > (1+\epsilon) \cdot \hat n_0(\lambda_5,\lambda_{6}) \quad & \mbox{and} \\
& \hat n_0(\lambda_j,\lambda_{j+1}) \leq (1+\epsilon) \cdot \hat n_0(\lambda_5,\lambda_{6}) \quad & \mbox{for all } i<j\leq4.
\end{eqnarray}
Then we put 
\begin{equation}
 \hat \beta_i = 1- \sum_{j=i+1}^6 \beta_j, 
\end{equation}
$\hat \beta_j=\beta_j$ for $j>i$ and $\hat \beta_1 = \ldots, \hat \beta_{i-1} =0$. \\
Case 2. If in case 1 there is no such index $i$ we may put $\hat \beta_j = \beta_j$ to be the pre-chosen weights for all 
$1\leq  j \leq 6$. 

The simulation of Section \ref{Section_Stationary_approach} is also conducted for the dynamic estimator (\ref{EsimatorDynamic}) 
with tuning parameter $\epsilon=0.05$ and the results are given in Table \ref{Table:001}. The corresponding dynamic adaptive 
SU test behaves almost like the adaptive SU test of the stationary approach, but shows significant advantages in the non-parametric 
case (D3). 

\begin{remark}
 The construction of the dynamic weights is motivated by (\ref{Pre-weighting}). Observe that the usual Storey estimator 
(\ref{StoreyEstimator}) mimics the ``slope'' of $t \mapsto \hat F_n(t)$ on the interval $[\lambda,1]$. The construction 
above splits that interval in smaller pieces with observed ``slope'' $\frac{\hat F_n(\lambda_j)- \hat F_n(\lambda_{j-1})}{\lambda_j - \lambda_{j-1}}$ 
on the interval $[\lambda_{j-1},\lambda_j]$. In comparison with the Storey estimator the backwards stopping procedure will 
stop when the observed slope increases too much. In conclusion the procedure may be viewed as a data driven selection of an 
estimation region $[\lambda_{i-1},1]$. Note that the new SU tests have finite sample FDR control. 
\end{remark}



% ...
%
%


\section{Adaptive SD tests}
\label{Section_Adaptive_SD_tests}

We now show that the results of the previous sections also hold for adaptive step down 
(SD) tests.  Again, let $\hat \alpha_{i:n}$, $i=1,\ldots,n$, be data dependent critical 
values with $0 = \hat \alpha_{0:n} < \hat \alpha_{1:n} \leq \ldots \leq \hat \alpha_{n:n} < 1$. Then the 
adaptive SD test rejects all null hypotheses corresponding to the $p$-values $p_{i:n}$ 
which fulfill $p_{i:n} \leq \hat\alpha_{R_{SD}:n}$ with 
\begin{equation}
 R_{SD} = \max\{i : p_{j:n} \leq \hat\alpha_{j:n} \mbox{ for all } j \leq i\} 
\end{equation}
and the convention $\max \emptyset =0$. Note that $R_{SD}$ again denotes the the number 
of rejections. Furthermore, let 
\begin{equation}
 V_{SD} = \sum_{i=1}^n 1\{p_i \leq \hat \alpha_{R_{SD}:n}, H_i=0\}
\end{equation}
be the number of falsely rejected true null hypotheses. 

\begin{theorem}\label{TheoremApaptiveSD}
 Consider the BI Model and the adaptive SD test with critical values (\ref{aSUcriticalvalues}) and 
estimator (\ref{aSUestimator}). Then 
\begin{equation}
 FDR = E\left( \frac{V_{SD}}{R_{SD} \vee 1} \right) \leq \frac{\alpha}{\lambda} E\left( \frac{V(\lambda)}{\hat n_0} \right) 
\end{equation}
holds. Moreover, the Theorems \ref{TheoremGeneralizedStoreyEstimator}, \ref{CorollaryStationaryApproach} 
and \ref{TheoremDynamicApproach} stay true if we replace adaptive SU test by adaptive SD test. 
\end{theorem}

The adaptive SU tests described in Section \ref{Section_FDR_control}-\ref{Section_Dynamic_approach} 
may also be carried out as adaptive SD tests with same critical values while still controlling the 
FDR at level $\alpha$. This has already be shown by \cite{Sarkar_2008} for the adaptive tests with 
estimators which may be treated by his FDR control condition. His methods include the Storey estimator 
(\ref{StoreyEstimator}), but exclude the generalized Storey estimators (\ref{GeneralizedStoreyEstimator}) 
with $\gamma_1<1$. 

\begin{remark}
 Since the distribution of false $p$-values is arbitrary, ties may occur. If $p_{j_1:n}=\ldots=p_{j_2:n}$ 
holds for some $j_1<j_2$, then for the computation of $R_{SD}$ these $p$-values are basically compared to 
$\hat\alpha_{j_1:n}$ which is the smallest of the corresponding critical values. However, Theorem 
\ref{TheoremApaptiveSD} also holds for the adaptive SD tests with data dependent critical values 
$\widetilde \alpha_{i:n} = \hat \alpha_{R(p_{i:n}):n}$, $i=1,\ldots,n$, where the $\hat \alpha_{i:n}$'s 
are given by (\ref{aSUcriticalvalues}) and $R(\cdot)$ is defined in (\ref{DefinitionRandV}). Then 
$p_{j_1:n}=\ldots=p_{j_2:n}$ are compared to $\alpha_{j_2:n}$ which is the largest of the corresponding 
critical values. The proof of this remark is included in the proof of Theorem \ref{TheoremApaptiveSD}. 
\end{remark}

\section{Converse Benjamini Hochberg Theorem}
\label{Section_Converse}

In this section we show under mild assumptions that the SU test of \cite{Benjamini_Hochberg_1995} is the 
only SU test with $FDR = \frac{E(N_0)}{n}\alpha$ under the BI Model. This converse BH-Theorem may be of 
interest when new testing procedures are supposed to be designed. The converse theorem relies on the 
following submodel of the BI Model given by Efron's mixture model (\ref{EfronsMixture}). 
\begin{itemize}
 \item[(C1)] The possible false $p$-values $\xi_1,\ldots,\xi_n$ are i.i.d. according to an unknown distribution 
	      function $F_1$ with Lebesgue density. In addition, let $\xi_1$ be stochastically smaller than a 
	      uniformly distributed $p$-value on the unit interval, i.e. $F_1(t) \geq t$, $0 \leq t \leq 1$. 
 \item[(C2)] The random variables $H_1,\ldots,H_n$ are Bernoulli variables and independent of $(\xi_i)_{i \leq n}$. 
	     Furthermore, let the range of probabilities $P(H_1=0)$ covers an interval $I\subset (0,1)$ with 
	     non-empty interior. 
\end{itemize}
Thus, the $p$-values $p_1,\ldots,p_n$ are i.i.d. according to the distribution function 
\begin{equation}\label{EfronsMixture}
 F(t) = P(H_1=0) \cdot t + P(H_1=1)\cdot F_1(t), \quad t \in [0,1]. 
\end{equation}

\begin{theorem}\label{TheoremConverse}
 Consider an adaptive SU test with critical values (\ref{aSUcriticalvalues}) and estimator (\ref{aSUestimator}). \\
(a) Suppose that $FDR = E(N_0) c$ holds for some constant $0<c<\frac{1}{n}$ for all distributions of the 
BI Model which satisfy (C1) and (C2). Then the sets of critical values have already the form 
\begin{equation}\label{TheoremConverse001}
\hat \alpha_{i:n} = \frac{i}{n}\alpha \mbox{ almost surely}, 
\end{equation}
for $i=1,\ldots,n$, where $\alpha = nc$. \\
(b) Suppose that $FDR = E(N_0) c$ holds for some constant $0<c<\frac{1}{n}$ for all distributions of the 
BI Model with fixed $N_0=1$ and (C1). Then (\ref{TheoremConverse001}) holds again. 
\end{theorem}

\begin{remark}
(a) There is no hope to obtain exact $FDR = \alpha$ for the adaptive SU tests for all distributions of the BI Model. 
The BH test is the only SU test with FDR which is ``distribution free'' with respect to the distribution of false 
$p$-values. \\
(b) In case of deterministic critical values $0<\alpha_{1:n} \leq \ldots \leq \alpha_{n:n} < 1$ we may choose 
$\lambda=\alpha_{n:n}$ and Theorem \ref{TheoremConverse} applies.  
\end{remark}

Suppose that $i \mapsto \frac{\alpha_{i:n}}{i}$ is non-decreasing. Then \cite{Benjamini_Yekutieli_2001} 
proved that the FDR is non-decreasing when the distribution $F_1$ of the false $p$-values becomes 
stochastically smaller (and the other way around when $\frac{\alpha_{i:n}}{i}$ does not increase). 
In these cases Theorem \ref{TheoremConverse} is not so surprising. However, the present result holds 
in general without any further monotonicity assumption of that kind and it is also true for our data 
dependent critical values.

%
%

\section{Proofs}
\label{Section_Appendix} 

{\sc Proof of Theorem \ref{TheoremCondition}.}
Conditioned under the generated $\sigma$-algebra 
$$\mathcal{F}_{\lambda} = \sigma(H_i,1\{p_i \leq s\} : s\geq \lambda, 1\leq i \leq n),
$$ 
the random variables $H_i$ and $1\{p_i\leq s\}$, $i=1,\ldots,n$, $s\geq t$, and in particular $n\hat F_n(s) = R(s)$, 
$s \geq \lambda$, can be treated as fixed deterministic values, due to measurability arguments. Under this 
condition we have exactly $V(\lambda)$ true $p$-values smaller or equal to $\lambda$, where $V(\lambda)$ is 
a fixed number too. Without restriction we assume $n_0 \geq V(\lambda)>0$ since everything is obviously fine 
for the excluded cases. Let us now consider new rescaled $p$-values $q_i$, $i=1,\ldots, R(\lambda)$ defined by 
$$ q_{i:R(\lambda)} := \frac{p_{i:n}}{\lambda}, \quad i=1,\ldots,R(\lambda).
$$
The new $p$-values $q_i$ corresponding to true null hypotheses are again i.i.d. and uniformly distributed on $(0,1)$ 
and independent from the rest of the $p$-values $q_i$ corresponding to false null hypotheses under the above condition. 
The exact positions of the $V(\lambda)$ true $p$-values in $(q_1,\ldots, q_{R(\lambda)})$ does not matter for our 
considerations. We now apply the BH Theorem for SU tests with critical values 
$$ \alpha_{i:R(\lambda)}^{(q)} := \frac{i}{R(\lambda)}\alpha' \quad 
\mbox{and } \alpha':= \frac{R(\lambda)}{\lambda \hat n_0} \alpha
$$
on the $q$'s. The data dependent level $\alpha'$ only depends by assumption (\ref{aSUestimator}) on the information 
given by $\mathcal{F}_{\lambda}$. Conditionally under $\mathcal{F}_{\lambda}$ we have a regular non data dependent 
SU procedure on the $q$'s. Let $R_q$ and $V_q$ denote the number of rejections and false rejections respectively 
by the above SU test. Observe that 
\begin{equation}\label{Theorem003Bew001}
 E\left( \frac{V_q}{R_q} \Big| \mathcal{F}_\lambda  \right) 
= \frac{V(\lambda)}{R(\lambda)} \min(\alpha',1) 
\end{equation}
holds by the BH Theorem. Obviously (\ref{Theorem003Bew001}) is $\frac{V(\lambda)}{R(\lambda)}$ in case $\alpha' \geq 1$.

Now observe that 
\begin{eqnarray*}
 R_q &=& \max\{i\leq R(\lambda) \, : \, q_{i:R(\lambda)} \leq \alpha_{i:R(\lambda)}^{(q)} \} \\
&=& \max\left\{i\leq R(\lambda) \, : \, \frac{p_{i:n}}{\lambda} \leq  \frac{i}{R(\lambda)} \frac{R(\lambda)}{\lambda \hat n_0} \alpha \right\} \\
&=& \max\left\{ i\leq n \, : \, p_{i:n}\leq  \left(\frac{i}{\hat n_0} \alpha\right) \wedge \lambda \right\} =R 
\end{eqnarray*}
and hence $V_q=V$ since both tests, belonging to $R$ and $R_q$, are rejecting the same hypotheses. 
Thus, by (\ref{Theorem003Bew001}) we get
\begin{eqnarray*}
    E\left( \frac{V}{R} \right) 
&=& E\left( E\left( \frac{V}{R}\, \Big| \mathcal{F}_\lambda  \right) \right)
= E\left( E\left( \frac{V_q}{R_q} \, \Big| \mathcal{F}_\lambda \right) \right)\\
&=& E\left( \frac{V(\lambda)}{R(\lambda)}\min(\alpha',1) \right) 
 =\frac{\alpha}{\lambda} E\left( V(\lambda) \min\left\{\frac{1}{\hat n_0}, \frac{\lambda}{R(\lambda) \alpha}\right\} \right).
\end{eqnarray*} 
$\hfill \square$

\par\bigskip
{\sc Proof of Theorem \ref{TheoremGeneralizedStoreyEstimator}.}
 First let us introduce the following simplifying notation 
\begin{eqnarray*}
& V(t)& := n_0 (t) := \# \{ p_i \leq t \, : \, p_i \mbox{ true}\} \\
& S(t)& := n_1 (t) := \# \{ p_i \leq t \, : \, p_i \mbox{ false}\},\quad 0\leq t \leq 1
\end{eqnarray*}
and similar to the proof of Theorem \ref{TheoremCondition} let us condition under 
$$\mathcal{F}_{\gamma_1} = \sigma(H_i,1\{p_i \leq s\} : s\geq \gamma_1, 1\leq i \leq n).
$$ 
Whereas $n_0(t)$ and $n_1(t)$ refer to known fixed values under the present condition, 
$V(t)$ and $S(t)$ are still random. Since $n_1(\gamma_1)-S(\lambda_1) \geq 0$ we obtain  
\begin{eqnarray}
&& \hspace{-0.75cm} E\left( \frac{V(\lambda)}{n(\hat F_n(\gamma_1) - \hat F_n(\lambda_1) +\frac{1}{n})} 
\Big| \mathcal{F}_\lambda \right)\\
&& \hspace{-0.75cm} = E\left( \frac{V(\lambda)}{n_0(\gamma_1)+1 -V(\lambda) - (V(\lambda_1) - V(\lambda)) + n_1(\gamma_1) - S(\lambda_1)
} \Big| \mathcal{F}_\lambda \right) \\
&& \hspace{-0.75cm} \leq E\left( \frac{V(\lambda)}{n_0(\gamma_1)+1 -V(\lambda) - (V(\lambda_1) - V(\lambda))
} \Big| \mathcal{F}_\lambda \right). \label{StarErsatz}
\end{eqnarray}
The random vector $(V(\lambda), V(\lambda_1)-V(\lambda), V(\gamma_1)-V(\lambda_1))$ is distributed according to 
the multinomial distribution 
$\mathcal{M}\left(n_0(\gamma_1), \frac{\lambda}{\gamma_1}, \frac{\lambda_1-\lambda}{\gamma_1},\frac{\gamma_1 - \lambda_1}{\gamma_1}\right)$ 
under our conditions. The subsequent Lemma yields
$$  (\ref{StarErsatz}) 
= \frac{\lambda}{\gamma_1-\lambda_1} \left( 1 - \left( \frac{\lambda_1}{\gamma_1} \right)^{n_0(\gamma_1)} \right)
\leq \frac{\lambda}{\gamma_1-\lambda_1}.
$$
Integration now gives 
$$  \quad \quad \quad \quad E\left( \frac{V(\lambda)}{n(\hat F_n(\gamma_1) - \hat F_n(\lambda_1) +\frac{1}{n})} \Big| (\hat F_n(t))_{t\geq \gamma_1} \right)
\leq \frac{\lambda}{\gamma_1-\lambda_1}. \quad \quad  \quad \quad \quad \square
$$

\par\bigskip
\begin{lemma}\label{LemmaAppendix001}
 Let $(V_1,V_2,V_3)$ be distributed according to the multinomial distribution $\mathcal{M}(n,p_1,p_2,p_3)$ with 
$p_3 >0$ and $n\geq 0$. Then we have 
\begin{equation}
 E\left( \frac{V_1}{n+1-V_1-V_2} \right) = \frac{p_1}{p_3} ( 1- (p_1 + p_2)^n).
\end{equation}
\end{lemma}

 \par\bigskip
{\sc Proof.}
 A simple calculation shows 
\begin{eqnarray*}
&& \hspace{-1cm} E\left( \frac{V_1}{n+1-V_1-V_2} \right) \\
&&= \sum_{{k_1> 0,\ k_2 \geq 0},\ {k_1+k_2 \leq n}} \frac{n!}{(k_1 - 1)! k_2! (n+1-k_1-k_2)!} p_1^{k_1} p_2^{k_2} p_3^{n-k_1 -k_2}\\
&&= \frac{p_1}{p_3} \cdot  \sum_{{j, k_2 \geq 0},\ {j+k_2 \leq n-1}} \frac{n!}{j! k_2! (n-j-k_2)!} p_1^{j} p_2^{k_2} p_3^{n-j -k_2}.
\end{eqnarray*}
The last equality follows from the substitution $j=k_1-1$. Observe that the last term adds the probabilities of the multinomial distribution 
times a constant factor, so extending the missing probabilities yields that the last term equals
\begin{eqnarray*}
\quad \quad && \hspace{-1cm} \frac{p_1}{p_3} \cdot \left( 1- \sum_{j+k_2=n} \frac{n!}{j! k_2! (n-j-k_2)!} p_1^{j} p_2^{k_2} p_3^{n-j -k_2} \right) \\
&& = \frac{p_1}{p_3} \cdot \left( 1- \sum_{j=0}^n \frac{n!}{j! (n-j)!} p_1^{j} p_2^{n-j} \right)
= \frac{p_1}{p_3} \cdot \left( 1- (p_1+p_2)^n \right). \quad \quad \square
\end{eqnarray*}

\par\bigskip
{\sc Proof of Lemma \ref{LemmaComparisonCondition}.}
Let us condition under $N_0=n_0>0$. Without restrictions assume that the $p$-values $p_1,\ldots, p_{n_0}$ belong to true 
null hypotheses. By (\ref{aSUestimator}) we have $\hat n_0^{(1)} = \hat n_0 $ on $\{ p_1 \leq \lambda \}$ and thus, 
\begin{eqnarray*}
&& \hspace{-1cm} E\left( \frac{V(\lambda)}{\hat n_0} \Big| N_0=n_0 \right) = \sum_{i=1}^{n_0}E\left( \frac{1\{p_i \leq \lambda \}}{\hat n_0} \Big| N_0=n_0 \right) \\
&& = n_0 E\left( \frac{1\{p_1 \leq \lambda \}}{\hat n_0} \Big| N_0=n_0 \right) 
= n_0 E\left( \frac{1\{p_1 \leq \lambda \}}{\hat n_0^{(1)}} \Big| N_0=n_0 \right) \\
&&= \lambda E\left( \frac{n_0}{\hat n_0^{(1)}} \Big| N_0=n_0 \right).
\end{eqnarray*}
The last equality holds because of the independence of $p_1$ and $\hat n_0^{(1)}$, where $p_1$ is uniformly 
distributed on $(0,1)$. There is nothing to show for $N_0=0$ and hence integration yields the assertion. $\hfill \square$

 \par\bigskip
{\sc Proof of Corollary \ref{CorollaryStationaryApproach}.}
The estimator $\hat n_0$ given in (\ref{WeightedStoreyEstimatorFixed}) is a convex combination. 
Thus, by Theorem \ref{TheoremGeneralizedStoreyEstimator} with $\gamma_i=1$ we obtain
$$ \quad \quad \quad E\left( \frac{V(\lambda)}{\sum_{i=1}^k \beta_i \hat n_{0,i}} \right) \leq 
\sum_{i=1}^k \beta_i  E\left( \frac{V(\lambda)}{\hat n_{0,i}} \right)
\leq \sum_{i=1}^k \beta_i \lambda = \lambda. \quad \quad \quad \quad \square
$$

 \par\bigskip
{\sc Proof of Theorem \ref{TheoremDynamicApproach}.}
Along the lines of the proof of the stationary approach (Corollary \ref{CorollaryStationaryApproach}) 
observe that 
$$ \frac{V(\lambda)}{\sum_{i=1}^k \hat\beta_i \widetilde n_0 (\gamma_i)} 
\leq \sum_{i=1}^k \hat\beta_i \frac{V(\lambda)}{ \widetilde n_0 (\gamma_i)} 
$$
holds for fixed $p$-values. Since $\hat\beta_i$ is $\mathcal{F}_{\gamma_i}$ measurable and 
(\ref{TheoremDynamicApproachCondition}), by taking conditional expectations we finally get 
\begin{eqnarray*}
 \quad E \left( \frac{V(\lambda)}{\hat n_0} \right) 
&\leq&  E \left( \sum_{i=1}^k \hat\beta_i \frac{V(\lambda)}{ \widetilde n_0 (\gamma_i)} \right) 
=E \left( \sum_{i=1}^k \hat\beta_i E \left( \frac{V(\lambda)}{ \widetilde n_0 (\gamma_i)} \Big| \mathcal{F}_{\gamma_i} \right)\right) \\
&\leq& E\left( \sum_{i=1}^k \hat\beta_i \lambda \right) = \lambda. \hspace{6cm} \square
\end{eqnarray*}

\par\bigskip
{\sc Proof of Theorem \ref{TheoremApaptiveSD}.}
(I) The proof is first done for the BH critical values $\alpha_{i:n} = \frac{i}{n} \alpha$ where also $\alpha \geq 1$ 
is allowed for SD tests at this point. For each $\alpha>0$ we will now prove 
\begin{equation}\label{LemmaSDProof001}
 E\left( \frac{V_{SD}}{R_{SD}} \right) \leq \frac{E(N)}{n}\alpha 
\end{equation}
for the BH SD test. Under the BI Model we can condition under $H_1,\ldots,H_n$ which are assumed to be fixed throughout. Then
\begin{equation}\label{Theorem1Proof001}
\begin{array}{ll}
 E\left( \frac{V_{SD}}{R_{SD}} \Big| H_1,\ldots,H_n \right) &= \sum_{i : H_i=0} E\left( \frac{1\{p_i \leq \alpha_{R_{SD}:n}\}}{R_{SD}} \Big|H_1,\ldots,H_n\right)\\
&= N_0 E\left( \frac{1\{p_1 \leq \alpha_{R_{SD}:n}\}}{R_{SD}} \Big| H_1,\ldots,H_n \right),
\end{array}
\end{equation}
where $p_1$ is assumed to be a true $p$-value, without restrictions. 
Let $p=(p_1,\ldots, p_n)$ be the vector of p-values and $p^{(1)}=(0,p_2,\ldots, p_n)$ the vector of p-values, where 
the first true $p$-value $p_1$ is decreased to zero.
Simple calculations show that $R_{SD}(p)=R_{SD}(p^{(1)})$ holds on the set $\{ p_1 \leq \alpha_{R_{SD}(p):n} \}$. 
Observe when $p_1$ is rejected and then $p_1$ could also be zero. Moreover we have 
$\{ p_1 \leq \alpha_{R_{SD}(p):n} \} \subset \{ p_1 \leq \alpha_{R_{SD}(p^{(1)}):n} \} $ in any case. Thus, 
\begin{eqnarray*}
 E\left( \frac{1\{p_1 \leq \alpha_{R_{SD}:n}\}}{R_{SD}} \Big| H_1,\ldots,H_n \right) 
&\leq& E\left( \frac{1\{p_1 \leq \alpha_{R_{SD}(p^{(1)}):n}\}}{R_{SD}(p^{(1)})} \Big| H_1,\ldots,H_n \right) \\
&\leq& E\left( \frac{\alpha_{R_{SD}(p^{(1)}):n}}{R_{SD}(p^{(1)})} \Big| H_1,\ldots,H_n \right) = \frac{\alpha}{n}
\end{eqnarray*}
by the independence of $p_1$ and $R_{SD}(p^{(1)})$ and Fubini's theorem. \\
(II) Inequality (\ref{LemmaSDProof001}) also holds for the modified data-dependent critical values 
$$ \widetilde \alpha_{i:n}' := \alpha_{R(p_{i:n}):n} = \hat F_n(p_{i:n}) \alpha.
$$
The same arguments of (I) can be applied by simply replacing $\alpha_{\,\cdot\,:n}$ by $\widetilde \alpha_{\,\cdot\,:n}$ and 
observing that $R_{SD}(p^{(1)}) = n\hat F_n(p_{R_{SD}(p^{(1)}):n})$ holds. 
\\
(III) As in the proof of Theorem \ref{TheoremCondition} we may condition under $\mathcal{F}_{\lambda}$. If we now use the 
inequality of (I) and (II) for the procedures, then Lemma \ref{TheoremApaptiveSD} follows from slightly adapted arguments given in 
the proof of Theorem \ref{TheoremCondition} and crossing over to inequalities. The arguments hold for $\hat \alpha_{i:n}$ 
as well as for $\widetilde \alpha_{i:n}$. The only difference for the SD case is that 
the upper bound of (\ref{Theorem003Bew001}) is $\frac{V(\lambda)}{R(\lambda)}\alpha'$ with the inequality ``$\leq$'' 
instead of ``='' in general. Furthermore $R_{q,SD} = R_{SD}$ 
still holds for both tests by an analogue argument for SD tests.
$\hfill \square$

\par\bigskip
{\sc Proof of Theorem \ref{TheoremConverse}.}
This proof uses results for complete statistical models which can be found in \cite{Lehmann_Romano_2005} 
and \cite{Pfanzagl_1994}. Our assumptions imply 
$$  \int \left( E\left( \frac{V}{R}\Big| N_0=n_0 \right) - c n_0 \right) \mathcal{L}(N_0)(d n_0) = 0
$$
for all distributions specified in (C1) and (C2). Observe that $\mathcal{L}(N_0)$ is a binomial distribution. 
Thus, 
$$  E\left( \frac{V}{R}\Big| N_0=n_0 \right) - c n_0 = 0
$$
holds for all $n_0 \in \{0,\ldots, n \}$ since $N_0$ is a complete statistic for the exponential family of 
binomials. Next we focus on $N_0=1$ and the proofs of (a) and (b) run parallel. Without restrictions we may 
assume that the true $p$-value is $p_1$ given by $U_1$. In that case $V$ and $R$ depend on the outcomes 
$(p_1,\ldots,p_n)=(U_1,\xi_2,\ldots,\xi_n)$, where $V(p_1,\ldots,p_n) = V(U_1, \xi_{1:n-1},\ldots,\xi_{n-1:n-1})$ 
holds and similarly for $R$. Note that here $\xi_{1:n-1},\ldots,\xi_{n-1:n-1}$ denote the order statistics 
of $\xi_2,\ldots,\xi_n$. Furthermore, put $\alpha=nc$. It is easy to see that $\alpha<1$ holds and the model 
includes i.i.d. uniformly on $(0,1)$ distributed $\xi_i$. Then we have 
$nE(\frac{V}{R}|N_0=1) = E(\frac{R}{R}1\{R>0\}) = P(R>0) < 1$ in both cases. Now we have 
\begin{eqnarray*}
 0 &=& E\left( \frac{V}{R} \Big| N_0=1 \right) - \frac{\alpha}{n} 
= E\left( \frac{V(U_1,\xi_2,\ldots,\xi_n)}{R(U_1,\xi_2,\ldots,\xi_n)} - \frac{\alpha}{n} \right) \\
&=& \int \int_0^1 \left( \frac{V(U_1, (\xi_{i:n-1})_{i\leq n-1})}{R(U_1, (\xi_{i:n-1})_{i\leq n-1})} - \frac{\alpha}{n}\right)\\
&& \hspace{4cm} 
d U_1 \ \mathcal{L}((\xi_{i:n-1})_{i\leq n-1}) d((\xi_{i:n-1})_{i\leq n-1}).
\end{eqnarray*}
Observe that the family of distributions of $\mathcal{L}(\xi_1)$ is convex and complete in the sense 
of \cite{Pfanzagl_1994} Theorem 1.5.10. Thus, the family of order statistics $(\xi_{i:n-1})_{i\leq n-1}$ 
is also complete and 
\begin{equation}\label{Theorem2proof1}
 \int_0^1 \frac{V(U_1, (\xi_{i:n-1})_{i\leq n-1})}{R(U_1, (\xi_{i:n-1})_{i\leq n-1})} \ dU_1 = \frac{\alpha}{n}
\end{equation}
 holds $\mathcal{L}((\xi_{i:n-1})_{i\leq n-1})$ almost surely. Consider first deterministic $\alpha_{i:n}$. 
Let $p=(p_1,\ldots, p_n)$ be the vector of p-values and $p^{(1)}=(0,p_2,\ldots, p_n)$ the vector of p-values, 
where the first true $p$-value $p_1$ is decreased to zero. Simple calculations show that $R(p)=R(p^{(1)})$ 
holds on the set $\{ p_1 \leq \alpha_{R(p):n} \}$. Observe when $p_1$ is rejected and then $p_1$ could also 
be zero. Moreover, we have $\{ p_1 \leq \alpha_{R(p):n} \} = \{ p_1 \leq \alpha_{R(p^{(1)}):n} \} $ in any 
case. The left hand side of (\ref{Theorem2proof1}) is thus equal to 
\begin{equation}\label{Theorem2proof2}
 \int_0^1 \frac{1\{p_1 \leq \alpha_{R:n}\}}{R} \ dp_1 = \frac{\alpha_{R(p^{(1)}):n}}{R(p^{(1)})}.
\end{equation}
The same arguments apply to the data driven critical values $\hat \alpha_{i:n}$. We have only to take 
$p_1 \leq \lambda$ into account. In this case we see that $\hat n_0 (p) = \hat n_0 (p^{(1)})$ does not 
change since $(\hat F_n(t))_{t\geq\lambda}$ is the same for $p$ and $p^{(1)}$. On $\{p_1\leq\lambda\}$ 
the $\hat \alpha_{i:n}$ can be considered to be deterministic and (\ref{Theorem2proof2}) also holds 
$\lambda \!\! \lambda_{(0,1)}^{n-1}$ $a.s$. It is easy to see that the set $\{R(p^{(1)})=j\}$ has 
positive probability for each value $1\leq j \leq n$ under at least one distribution of model (C1). 
Therefore, consider uniformly on $(0,1)$ distributed false $p$-values and observe that 
$$ \{ p_{j:n} \leq \alpha_{j:n}, \, p_{j+1:n} > \alpha_{n:n} \} \subset \{R(p^{(1)})=j\} 
$$
holds for deterministic critical values and 
$$ \{ p_{j:n} \leq \hat\alpha_{j:n}, \, p_{j+1:n} > \lambda \} \subset \{R(p^{(1)})=j\} 
$$
for our data driven critical values. 
From (\ref{Theorem2proof1}) and (\ref{Theorem2proof2}) we conclude the equality $
 \frac{\alpha_{j:n}}{j} = \frac{\alpha}{n}$ for all $1\leq j \leq n.
$ $\hfill \square$




\section*{Acknowledgements} 

The authors are grateful to Helmut Finner who introduced us in the field of multiple testing. We 
also wish to thank Julia Benditkis for helpful discussions and the Deutsche Forschungsgemeinschaft 
(DFG) for financial support.

 \bibliographystyle{elsarticle-harv} 
 \bibliography{LiteraturFDR.bib,LiteraturSonstiges.bib}


\end{document}